\documentclass{amsproc}

\usepackage{graphicx}
\usepackage{xcolor}

\definecolor{darkgreen}{rgb}{0,.5,0}
\definecolor{brown}{rgb}{0.5,0.3,0}

\usepackage{silence}
\usepackage[hyphens]{url}
\usepackage{url}

\theoremstyle{definition}

\theoremstyle{remark}

\numberwithin{equation}{section}

\begin{document}

\title[Todxs cuentan: community and belonging from the first day of class]{Todxs cuentan:  building community and \\ welcoming humanity from the first day of class.}

\author{Federico Ardila--Mantilla}
\address{San Francisco State University and Universidad de Los Andes}
\email{federico@sfsu.edu}

%\subjclass[2010]{Primary } I'll ask for these later if we need them.
%\keywords{}

\date{\today}

%\dedicatory{Optional dedication}%optional

\begin{abstract} 
Everyone can have joyful, meaningful, and empowering academic experiences; but no single academic experience is joyful, meaningful, and empowering to everyone. How do we build academic spaces where every participant can thrive? 
Audre Lorde advises us to use our differences to our advantage. bell hooks highlights the key role of building community while addressing power dynamics. Rochelle Guti\'errez emphasizes the importance of welcoming students' full humanity. 
This note discusses some efforts to implement these ideas in a university classroom, focusing on the first day of class.
\end{abstract}

\maketitle

\section{Community.}

\begin{quote}
Excitement about ideas [is] not sufficient to create an exciting learning process.
As a classroom community, our capacity to generate excitement is deeply affected by our interest in one another, in hearing one another's voices, in recognizing one another's presence. Any radical pedagogy must insist that everyone's presence is acknowledged. That insistence cannot be simply stated. It has to be demonstrated through pedagogical practices. There must be an ongoing recognition that everyone influences the classroom dynamic, that everyone contributes. Often before this process can begin there has to be some deconstruction of the traditional notion that only the professor is responsible for classroom dynamics.

\begin{flushright}
bell hooks  \cite{hooks}
\end{flushright}

\end{quote}

\section{January, 2017: Before classes start.}

The week before a semester starts, I often find myself frantically trying to organize my office, our apartment, our record collection, anything else that needs or does not need organizing. This makes me feel productive while I avoid preparing for my upcoming classes. One day, organizing our living room, I found a portable turntable and far too many records that did not fit in our crates. I brought them to my office. When I played the first song, I was pleasantly surprised by how drastically this addition transformed the space -- just like the extra coffee maker I had found a few semesters ago, which now allows me to offer visitors a nice, strong cafecito before we begin to talk about life, or mathematics, or both.

The first days of class heavily dictate how a classroom will feel throughout the semester. The physical setting also plays a crucial role. I knew this semester I would be teaching in a dark room with small windows and broken blinds. The whiteboards on all walls would be very useful for group work; but the long rows of tables nailed to the ground, and the clunky laptop computer locked into place on every seat, would make collaboration challenging. 

A few hours before class, I was thinking about how to make students feel welcome in this space. This was a combinatorics class with a very broad range of students: from second-year undergraduates to Master's students doing research in the field. Most of them did not know each other, and I was dreading the uncomfortable silence that can sometimes engulf the room before class starts. So I thought: ``I should at least bring my turntable and a few records". 

When students arrived on the first day of class, Carlos Embales was playing. This quickly broke the ice, and seemed to give them permission to start talking to each other.

\begin{figure}[h]
\includegraphics[height=6cm]{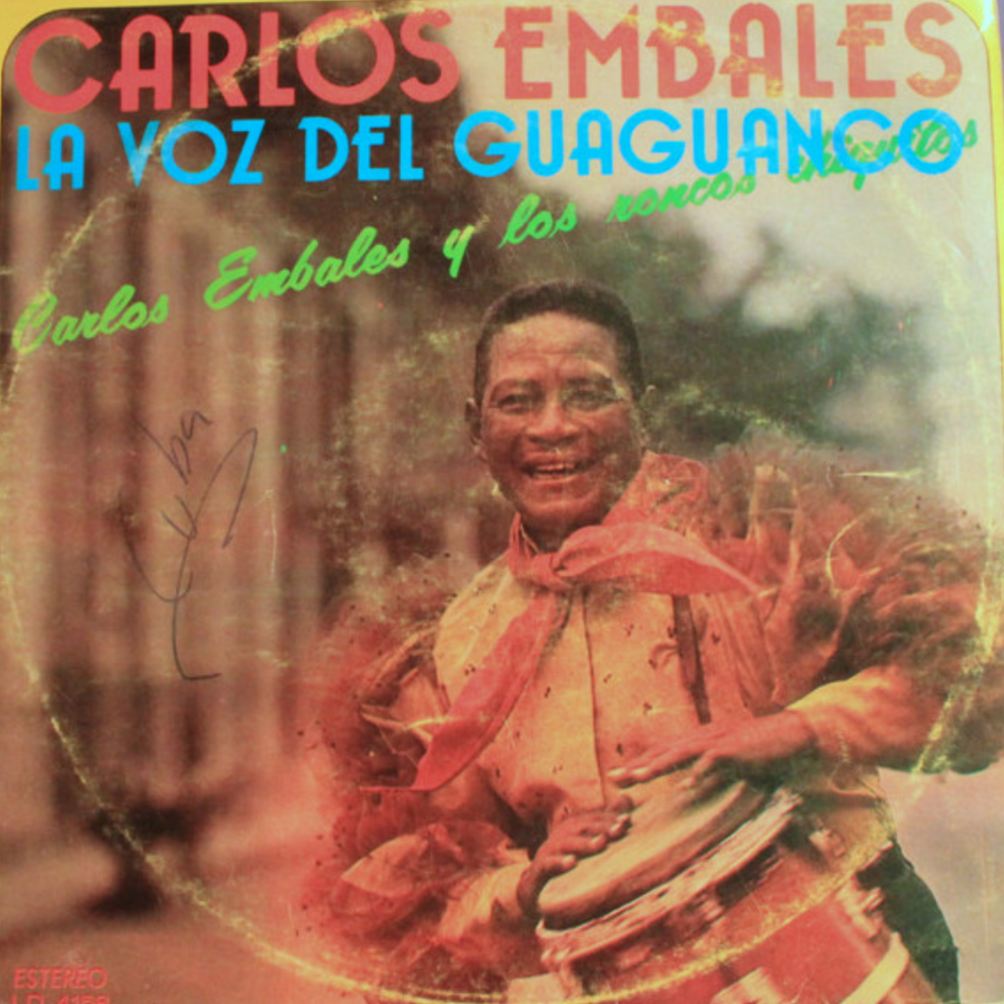}
\end{figure}

\section{Introductions.}

We got in a circle 
and answered a few questions: 

\noindent $\bullet$
What would you like us to call you? 

\noindent $\bullet$
What is something outside of mathematics that you love doing? 

\noindent $\bullet$
How do you feel about being here?

I answered first. I like being called Federico or Profe. I DJ. I am excited and a bit nervous, because I am going to try many new things in class this semester, and although academic tradition dictates that a professor is supposed to appear invulnerable and in control, I plan to put us in learning experiences that will not really be under my control.

My students love making music, dancing, designing, playing video games, trying to solve  crimes where the suspect was wrongly convicted. They mostly said they were excited and nervous, like me.

I explained my vision for this exercise: we say these things out loud to remember that mathematics is a human endeavor. I wanted to make clear that our full humanity was not only welcome here, but in fact would define and enrich this mathematical space.

\section{Qu\'{\i}talo del Rinc\'on.} 

After these introductions, I played a song:
\begin{center}
\url{http://math.sfsu.edu/federico/Talks/embale.m4a}
\end{center}

\noindent In case you are not able to hear it, let me try to describe it, certainly not doing it any justice. The song starts with a pair of sticks, three drums, and a shaker, weaving an intricate combination of rhythms. A singer chants a long, melodic \emph{aaalalalalalalalaaa}. Then he is joined by the main singer; in a beautiful and mysterious harmony, they introduce the theme of the song. The chorus comes in: a joyful call and response between the lead singer and a group of high-pitched voices -- kids, maybe. While the kids keep repeating the chorus, the main singer starts improvising rhymes, and the drums take off. If you have been a part of a group like this, you will recognize the feeling: you stop knowing exactly what it is that you are doing, and you collectively connect to something deeper than anything you can reach on your own. After a couple of minutes the recording fades out, but you can tell this is just for technical or commercial reasons: the musicians show no sign of slowing down; they are only getting started.\footnote{A recent performance of this song is at \url{https://www.youtube.com/watch?v=FpxE_xzPNQY.}}

As the music played, I asked students to come up with a mental picture of what was happening, and write down a few words to describe it. I'll invite you to do that as well. Because students who could understand the lyrics were better prepared to answer this question, I asked them to step back for a moment and let the others answer first. They said:

\medskip

\begin{quote}
\begin{center}
 community . joy . polyrhythm . family . crescendo . playful \\
 encouraging . unexpected . churchlike . inviting . dancing  \\
 conversation . courage . motivation . cheerful . Spanish \\
 learning . rhythm . celebration . style . culture . festive 
\end{center}
\end{quote}

\medskip

Eventually, the Spanish speakers in the class told the rest what is happening here: It is a math class!
 \emph{Qu\'{\i}talo del rinc\'on} by Carlos Embales y los Roncos Chiquitos is a guaguanc\'o -- the style of Cuban rumba native to the Black neighborhoods of La Habana, born soon after the abolition of slavery in the late 1800s. The chorus says:

\medskip

\begin{quote}
\begin{center}
If someone doesn't wanna learn, we'll teach them, very happily! 

\smallskip

20+3? \textbf{23}. \quad 30+6? \textbf{36}. \quad
20+3? \textbf{23}. \quad 30+6? \textbf{36}. \quad 

\smallskip

Take them away from the corner, bring them to the window;

you'll see how they'll learn right away, full of joy!

\smallskip

20+3? \textbf{23}. \quad 30+6? \textbf{36}. \quad
20+3? \textbf{23}. \quad 30+6? \textbf{36}. \quad 
\end{center}
\end{quote}

\bigskip

On break during a mathematical visit to Ann Arbor, MI, I found this enigmatic album for \$1, digging through the sales bin of a music store -- the dream of a record collector! I was even more thrilled when I got back home and I heard this song. I return to it often, when I think about what I'd like my math classroom to feel like.

\section{A Community Agreement}

The course syllabus is the first official document students receive in a class; it is the first impression they receive about what is valued in the class. When I began teaching, I would make the syllabus the night before classes started, essentially copying the syllabus from whoever taught the class last time -- including the grading scheme. As an unintentional consequence, my class often valued whatever the last instructor valued. 

In recent years I have tried to write syllabi that actually communicate the kind of course that I hope to build together with my students. This semester, after playing and discussing \emph{Qu\'{\i}talo del Rinc\'on}, we discussed the first part of the syllabus, which read:

\begin{figure}[h]
\includegraphics[height=4.5cm]{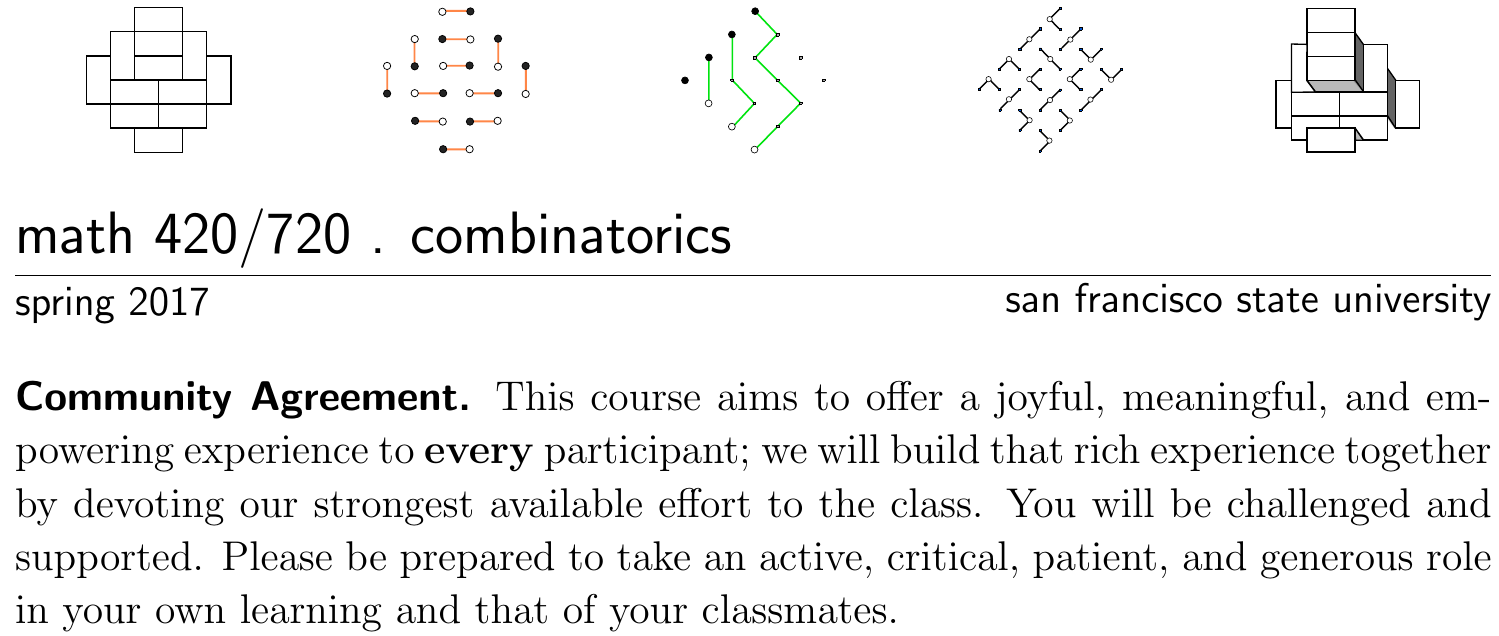}
\end{figure}

In small groups, students discussed this agreement. To initiate a dialogue about it, I asked each group to share two words that particularly resonated with them. The wide variety of different answers was striking to me. Some students were excited that they would be challenged; some that they would be supported; some liked the combination of the two. Many liked the word ``available"; most of them work and have many other obligations, and they appreciated that this was acknowledged. We discussed how to be productively critical of each other's work, and what generosity might mean in a mathematics classroom. We talked about how sometimes we are very good at being patient with our classmates, but we are not so good at being patient with ourselves.

My students and I have cocreated this Community Agreement over the last few years. Each semester, students start with the one used in the previous semester, and they make some (usually small, always thoughtful) changes. I then incorporate these changes into future agreements. 

I think it is important that this does not feel like an externally imposed code of conduct that they must obey. Instead, my hope is that we can reach a community agreement that is actually ours, that we are all excited to put into action.

\section{Assessment}

To conclude the first day of class, we discussed the grading scheme for the class. 
Typically, grades in math classes have been based on exams and, to a lesser extent, on homework. This disproportionately rewards a certain kind of mathematician, who enjoys and thrives solving problems quickly under pressure and time constraints. I am that kind of mathematician, and I suspect that so are most current professors of mathematics. After all, we had to succeed in this kind of grading system to become professors. But if we are honest with ourselves -- and I am honest with my students about this -- this is a very narrow kind of assessment; one that is easy for us to grade, but is not actually very good at measuring the kind of deep, creative thought that we associate with mathematical work.

\begin{quote}
We value what we measure because we do not know how to measure what we really value.
\begin{flushright}
Richard Tapia \cite{Tapia}
\end{flushright}
\end{quote}

Over the years, I have asked students for suggestions on the kinds of assessment that will most accurately reflect the mathematical work they are capable of. I have since shifted to grading schemes that seek to reward many types of mathematical work. For this paired undergraduate/graduate course, it consisted of: individual and group homework, daily notecards (in pairs) and a ``math diary" (individually) summarizing the main lessons learned each day, optional investigations on open problems, and a small research project (preferably but not mandatorily in groups). There were no exams; I have always found combinatorial reasoning especially difficult to come up with on the spot.

\section{From Abstract Goals to Concrete Practices}

On Day 1, my students and I generally concurred that our Community Agreement, and the words they used to describe \emph{Qu\'{\i}talo del Rinc\'on}, were good goals to aim for in our class. They felt that

\begin{quote}
\begin{center}
 community . joy . polyrhythm
\end{center}
\end{quote}

\noindent were especially important. 
With that in mind, I asked them to come up with a few concrete practices that we could follow to build this classroom culture and atmosphere together. They wrote them in notecards and brought them to class on Day 2, anonymously if they preferred.

\medskip

They had some suggestions for me, the instructor; for example:

\begin{itemize}
\item
Offer many group assignments where we get to work with different people each time.
\smallskip
\item
You told us that the course will emphasize growth and teamwork. Have the assessment and the grading reflect this. 
\end{itemize}

\medskip

Students also had many suggestions for themselves and each other:

\begin{itemize}
\item
Let's be very mindful of how we communicate with each other.
Emphasize constantly that mathematics is often difficult, and understanding is developed through extensive practice. Replace ``this is obvious" with ``with a bit of thought one can understand this"; ``I'm stupid" with ``I'm struggling"; ``I can't do this" with ``I can't do this yet". 
\smallskip
\item
Let's not take the joy of discovery away from others. If I think I understand something, I should step back for a moment, and offer myself as a resource to others as needed. 
\smallskip
\item
Let's stay honest and vulnerable. If I don't think I understand something, I should ask for help. 
\smallskip
\item
Let's be excited to help our classmates learn, with some leadership from the teacher.
\smallskip
\item
It was so interesting how every instrument plays a totally different rhythm but altogether they create a very beautiful piece of art. Similarly, every brain works differently, and creating a math community to solve problems will make learning much fun, and will lead to more creativity.
\smallskip
\item
In that \emph{guaguanc\'o} we can only hear the musicians, but we're pretty sure the community is dancing right in front of them. Try to accomplish that in our class.\footnote{I still think a lot about how we might do this.}
\end{itemize}

\smallskip

The concrete suggestion that we most often returned to was the following:

\begin{itemize}
\item
\textbf{Make space, take space.} If I feel comfortable speaking out, I should be mindful of how much space I take, and make room for others. If I tend to be quieter in groups, I should remember that my ideas are important, and others will benefit from hearing them.
\end{itemize}

\smallskip

All of these suggestions became part of our course syllabus.

\section{Make space, take space.} 

This last piece of advice feels really relevant to me as I write this. 

As a mathematics researcher with more than 20 years of experience, I feel pretty confident that my mathematical ideas are valuable. It sometimes takes a special effort to truly listen to students' ideas without projecting my own views onto them. When I have been able to really make space for students' thought, we have all learned very innovative and useful ways of thinking about combinatorics. 

As a mathematics educator with great interest but under 20 minutes of formal training in education, I still feel like a student with everything to learn. Writing about pedagogy feels very uncomfortable. For every criticism the reader may have of my educational work, I have at least five. I cannot count the number of self-deprecating statements I have edited out of this note.

However, I did commit to upholding our Community Agreement. Now that the editors of the volume have made space for my thoughts, I feel compelled to embrace our collective cultural practice, take space, and speak-while-uncomfortable anyway.

\section{Music}

After our first meeting it occurred to me that, if I was asking students to help me create an ideal atmosphere for our class, then I should not be the only one choosing music for us. Thus the first homework read:

\begin{quote}
Homework 0. Let's continue playing some music before class, to bring some more light into that classroom. On your designated day, please choose a song to share that makes you feel comfortable, joyful, at home. If you'd like to, you can tell us a bit about the song or why it's meaningful to you.
\end{quote}

\noindent We did this throughout the semester.

A. got us started, playing a live performance of \emph{ya ba7riye/shiddo el-himmi}
by Marcel Khalife; in the chorus, a stadium full of people sings ``Oh freely, hey hey hey hey".
She told us that as a Palestinian woman and an immigrant in the US, raising four children while working and going to school at the same time, it's very difficult for her to feel at home and welcome and free in this country, in this day and age. But mathematics is a place where she feels free, where no one can take her freedom away.

R. gave the class a three-song showcase of cumbia's migration from Colombia to Mexico to California. ``Every Californian should know about cumbia."

T. chose \emph{Dear Mama} by 2Pac: ``My mom worked incredibly hard to give me the opportunity to go to college; when I'm in these classrooms, I am constantly thankful to her." This clearly resonated with several students, singing along.

C., a software engineer turned mathematician, shared the music and the journey of software engineer turned singer-songwriter Vienna Teng.

B. and C. and J. shared their favorite songs to perform. 

C. told us that she wanted to share the song she sang at her mother's funeral. We did our best to hold space for her.

J. played Lauryn Hill; ``Who doesn't want to hear Lauryn Hill?"

M. made sure we knew that the Filipino-American hip-hop scene in the Bay Area is still going strong. Some of us knew Rocky Rivera as an MC, but none of us realized that she was a student on our campus.

A. brought her daughter M. to class, and played her favorite song: Israel `IZ' Kamakawiwo`ole's \emph{Somewhere Over The Rainbow}. It assures young M. that all her dreams can come true.

I was sure that my students would bring lots of good music, but I would have never imagined how deeply personal this exercise turned out to be. One thing seemed clear to me: my students wanted to be seen, really seen, as full humans, inside the classroom.

\section{How I experienced this classroom.}

Throughout my career I have tried to make my mathematics classroom a human place, where every interested student feels at home, and finds a conducive environment to discover and shape their own mathematical voice. It's a tall order, and I certainly will not claim that I have succeeded.

I will say this: That SFSU Combinatorics class felt like no other that I'd ever experienced. Teaching and learning it was a tremendously human experience for me. Additionally, and relatedly, this was also the home to the richest mathematical discussions I had ever seen in one of my classes.

Let me confess something.  
When we spend a whole hour getting to know each other, 
when we spend a few minutes of every class sharing music that is meaningful to us, 
when we spend most of the class hour exploring mathematical situations together and at most 15 minutes ``delivering content", 
I start worrying: Am I covering enough mathematics? 

I have come to understand that when students are engaged so actively, and when we really listen to each other's ideas, a creative, mathematical magic can happen that I could not have arrived at by simply preparing a lecture and delivering it. In this class, more than ever before, I experienced my students truly take charge of their shared learning experience, take ownership of the material, allow themselves to ask their own critical, insightful mathematical questions, value those questions, and turn them into their own original discoveries. In fact, their insight taught me many new things about classic problems that I thought I understood completely. More importantly, it led to new discoveries that I think only they would have come up with.\footnote{The mathematical work of this classroom will be the subject of an upcoming paper.}

I cannot take credit for this. In fact, I am certain that I will not be able to replicate it: a unique combination of humans made this classroom what it was, and led to a unique atmosphere and a unique mathematics. In the future, as a professor, I can only try to put some structure in place that may help my students and I flourish together. I continue to do this, each semester, with varying success.

\section{How students experienced this classroom.}

\begin{quote}
Teachers cannot claim their pedagogy is rehumanizing without obtaining recurring evidence from their students that they agree and without giving students opportunities to offer additional approaches for rehumanizing.
\begin{flushright}
Rochelle Guti\'errez \cite{Gutierrez}
\end{flushright}
\end{quote}

The following is a representative selection of students' feedback on the course.

\begin{itemize}

\item
The first day of a class wasn't spent reading through the syllabus or diving into material. Rather, it was spent entirely on introductions and conversation, setting the tone for a class in which students are deeply valued as human beings rather than just as mathematicians.

\smallskip
\item
The math was great, but the thing that stood out to me was the music. As I have been teaching for now 4 years, I try to continuously find different ways to make students feel comfortable/motivated/etc.  [...] having everyone have a chance to express themselves in that way in the class was awesome, so awesome that I actually used it in my class this semester.
 
\smallskip
\item
I am totally stealing classroom structures used this semester to
implement in the classes that I will teach in the future.  
 
\smallskip
\item
These are the kind of classes that remind me why I love math. I really enjoy the learning environment that was created. [...] I am not a fan of group work with preassigned groups. I would keep getting in the same group with an individual who's learning style was less than compatible with my own which was frustrating and a bit unnecessary.

\smallskip
\item
He tries really hard to engage with everyone and that paradoxically means that he doesn't have a lot of time for an individual student sometimes.
 
\item
I'm typically one who doesn't speak out much in class but working in groups helped me to become more comfortable and I found myself sharing more than usual.
 
\smallskip
\item
He ensures that each class member knows their opinions are important and that their voices should be heard. We established a supportive atmosphere and frequently worked in groups on difficult and interesting problems making sure everyone made a significant contribution and had a strong grasp on the material. I found myself pouring my extra time into this class because of how much I enjoyed learning the material we covered.

\smallskip
\item
As a combinatorics enthusiast, I have seen or self-discovered all problems/techniques covered in the class. I cannot afford to spend my sharpest years not learning, especially if I want to contribute to combinatorics. Each day, I hoped for something new, but each day was disappointing. I really did enjoy the homework, but I stopped attending class. Instead, I read the book at a higher rate at home and self studied upon completing the book.

\smallskip
\item
Math departments can be inhospitable elitist places where undergraduates who are earlier in their careers are looked down upon for not immediately grasping concepts.
From the first day of class he builds a supportive environment for those students who may feel "non-brilliant" and helps them see that they have just as much to contribute as other students. He is always willing to seriously consider a student's ideas and suggestions.
All of this while still pushing each of us to challenge ourselves and providing ways for students to pursue their specific interests. 

\end{itemize}

\section{Acknowledgments.}

This work was partially supported by NSF grant DMS-1855610, a Simons Fellowship, and NIH SF BUILD grant 5UL1GM118985-03. I began writing it during an SF BUILD Writing Retreat in 2017, and finished it while on sabbatical at the Universidad de Los Andes in Bogot\'a, Colombia in 2020. I thank SF BUILD, SFSU, Los Andes, and the Simons Foundation for providing excellent working conditions. In particular, I wish to thank Camille Rey for her excellent writing advice.

I am very grateful to the numerous people who have shaped my ideas on mathematics education, and to those who have encouraged me to share them. In particular, I am very grateful for the invitation to contribute to this volume. 

My pedagogical practices owe tremendously to my colleagues at SFSU, particularly Kim Seashore and Kimberly Tanner.
I have learned a lot about supporting students of color in science from the SF BUILD Faculty Agents of Chance program. 
Audre Lorde, bell hooks, and Rochelle Guti\'errez have shaped my views on education. 
Carlos Embales gave me a priceless and unexpected pedagogy lesson. 

Most importantly, I would like to thank my wonderful students; they are my teachers. Their diverse, critical, and generous perspectives have completely transformed my world view and my understanding of our work as educators. Working with them keeps my spirit young and my heart full.

\end{document}